\documentclass{article}
\usepackage{epsfig}
\usepackage{latexsym}
\usepackage{amsmath}
\usepackage{amsthm}
\usepackage{amssymb}
\usepackage{graphicx}
\usepackage{psfrag}
\usepackage{thesis}
\usepackage{hodge}
\usepackage{eufrak}
\begin{document}
\title{Generating Functions for Hurwitz-Hodge Integrals}
\author{Renzo Cavalieri}

\maketitle
\begin{abstract}
In this paper we describe explicit generating functions for a large class of Hurwitz-Hodge integrals. These are integrals of tautological classes on moduli spaces of admissible covers, a (stackily) smooth compactification of the Hurwitz schemes.

Admissible covers and their tautological classes are interesting mathematical objects on their own, but recently they have proved to be a useful tool for the study of the tautological ring of the moduli space of curves, and the orbifold Gromov-Witten theory of DM stacks.
 
Our main tool is Atiyah-Bott localization: its underlying philosophy is  to translate an interesting geometric problem into a purely combinatorial one. 
\end{abstract}
\section*{Introduction}

\subsection{Informal Overview}
For any fixed integer $d$ consider
$$\mathcal{A}_g:=\A{g}{0}$$
the moduli space of genus $g$ admissible covers of a rational curve, with two fully ramified points (see Section \ref{adm}). $\mathcal{A}_g$ is a smooth Deligne-Mumford stack of dimension $2g-1$. The goal of this paper is to evaluate, for all pairs of integers $(j,g)$, the tautological classes (see Section \ref{taut}):
$$
\lambda_g \lambda_j \psi_{2g+2}^{g-j-1}.
$$
\textbf{Notation (abuses of):} 
\begin{enumerate}
	\item Since we are always dealing with $\psi$ classes at the (fully ramified!) last point, we immediately drop the subscript $2g+2$ from our notation. 
	\item To try and keep this overview more legible, we  omit the integral signs in this section. Every time we talk about a dimension $0$ tautological class, we implicitly mean its evaluation on the appropriate moduli space.
\end{enumerate}
 
Our intention is to describe generating functions for such numbers; therefore, we attach the formal variable $u^{2g}$ to any integral computed on the space $\mathcal{A}_g$.

\begin{figure}
	\begin{center}
	\psfrag{D1}{$\mathcal{D}_1$}
	\psfrag{D2}{$\mathcal{D}_2$}
	\psfrag{D3}{$\mathcal{D}_3$}
	\psfrag{D4}{$\mathcal{D}_4$}
	\psfrag{T1}{$\mathcal{T}_1$}
	\psfrag{T2}{$\mathcal{T}_2$}
	\psfrag{T3}{$\mathcal{T}_3$}
	\psfrag{T4}{$\mathcal{T}_4$}
	\psfrag{V0}{$\mathcal{V}_0$}
	\psfrag{V1}{$\mathcal{V}_1$}
	\psfrag{V2}{$\mathcal{V}_2$}
	\psfrag{V3}{$\mathcal{V}_3$}
	\psfrag{u0}{$1$}
	\psfrag{u2}{$u^2$}
	\psfrag{u4}{$u^4$}
	\psfrag{u6}{$u^6$}
	\psfrag{u8}{$u^8$}
	\psfrag{0}{$\frac{1}{d}$}
	\psfrag{1}{$\lambda_1 $}
	\psfrag{2}{$\lambda_2 \psi$}
	\psfrag{3}{$\lambda_3 \psi^2$}
	\psfrag{4}{$\lambda_4\psi^3$}
	\psfrag{21}{$\lambda_2\lambda_1$}
	\psfrag{31}{$\lambda_3\lambda_1\psi$}
	\psfrag{41}{\hspace{-0.03cm}\small{$\lambda_4\lambda_1\psi^2$}}
	\psfrag{32}{$\lambda_3\lambda_2$}
	\psfrag{42}{$\lambda_4\lambda_2\psi$}
	\psfrag{43}{$\lambda_4\lambda_3$}
		\includegraphics[width=.95\textwidth]{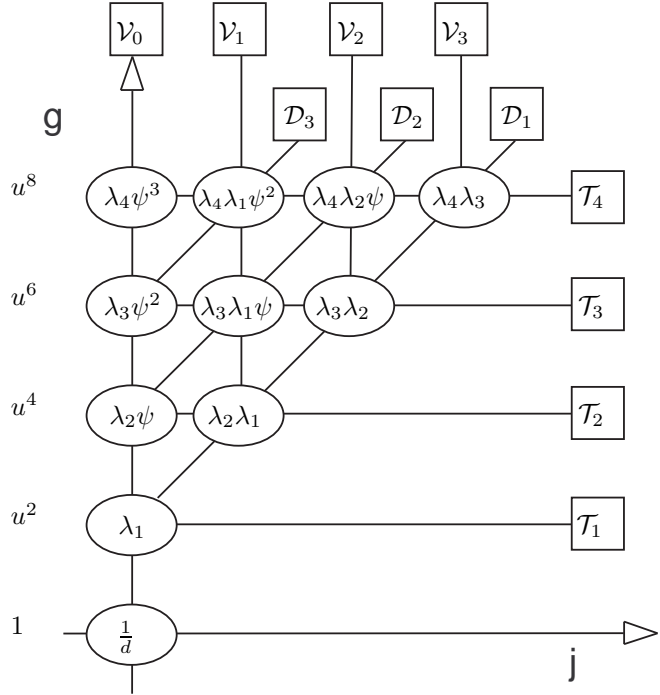}
	\end{center}
	\caption{A sketch of our generating function convenctions.}
	\label{dd2}
\end{figure}

In Figure \ref{dd2}, we  sketch the possibly non zero numbers on the $(j,g)$-plane. We highlight three natural ways of grouping our numbers into generating functions:
\begin{description}
	\item[Diagonally:]  add all terms that lie on lines parallel to the diagonal $j=g$. We define $\mathcal{D}_i(u)$ to be the generating function corresponding to the line $j=g-i$.
	\item[Vertically:]  fix $j$ and vary $g$. We denote by $\mathcal{V}_i(u)$ the function obtained by fixing $j=i$.
	\item[Horizontally:] fix $g$ and vary $j$. In this case the power of the formal variable remains the same. We define $\mathcal{T}_g(u)$\footnote{unfortunately $\mathcal{H}$ is already taken by ``Hurwitz". We choose $\mathcal{T}$ as in ``Total Chern class", since $\mathcal{T}_g(u)$ can be viewed as $u^{2g}\lambda_g c(\mathbb{E})/(1-\psi)$.} to be $u^{2g}\lambda_g\sum_j \lambda_j \psi^{g-j-1}$.
\end{description} 
We now conclude this section by paraphrasing the main results of this paper, referring the reader to Section \ref{results} for the precise statements.
\begin{description}
 	\item[Paraphrase of Theorem \ref{d}:]\textit{ The generating function $\mathcal{D}_i$ can be obtained from $\mathcal{D}_1$. In particular $\mathcal{D}_i$ is the $i$-th term in the expansion in $\mathcal{D}_1$ of $\frac{1}{d}$e$^{d\mathcal{D}_1}$.}
 	\item[Paraphrase of Corollary:]\textit{The function $\mathcal{T}_g$ is the degree $2g$ term in the series expansion of
 	$$
 	\frac{d^{d-1}\sin^d\left(\frac{u}{2}\right)}{\sin^d\left(\frac{du}{2}\right)}.
 	$$}
 	\item[Paraphrase of Theorem \ref{v}:]\textit{ The generating function $\mathcal{V}_i$ can be reconstructed from the values of $\lambda_g\lambda_{g-1}$, for $1\leq g\leq i$. The formula is given explicitly, and it is a product of an exponential function with a polynomial indexed by partitions of the integer $i$.}
\end{description}
\subsection{History, Connections and Applications}
\label{history}

Hodge integrals are evaluations of certain dimension zero tautological classes on $\overline{M}_{g,n}$: polynomials in $\lambda$ classes and $\psi$ classes. Besides being interesting mathematical objects on their own, their importance lies in the fact that they seem to create a powerful connection between three areas of mathematics: the study of the geometry of the moduli space of curves, the combinatorial/representation theoretic Hurwitz theory, and the Gromov-Witten theory of toric varieties.

The ELSV formula (\cite{elsv:hnaiomsoc} and \cite{elsv:ohnahi}) states that Hurwitz numbers can be expressed as Hodge integrals. This has been the springboard for interesting work of Tom Graber, Ravi Vakil, and the combinatorialists Ian Goulden and David Jackson (\cite{gv:taut}, \cite{gv:hnavl}, \cite{gv:rvl}, \cite{gjv:ttgodhn}, \cite{gjv:spolgc}), making progress towards the understanding of the celebrated Faber conjectures on the tautological ring of $M_{g,n}$ and various (partial) compactifications of it.

Atiyah-Bott localization provides the link with Gromov-Witten theory: when localizing on the spaces of stable maps to a toric variety, the fixed loci are essentially\footnote{Often one has to consider also quotients by finite groups, brought in by the presence of automorphisms of maps.} products of moduli spaces of curves. The restriction of the virtual fundamental class to the fixed loci and the contribution of the euler class of the normal bundle are expressed in term of $\lambda$ and $\psi$ classes, thus giving rise to products of Hodge integrals. 

In \cite{f:algo}, Carel Faber provides an algorithm for computing any given Hodge integral. However, we know only of very few examples of classes of Hodge integrals that are explicitly described in generating function form. This is again work of Faber and Pandharipande (\cite{fp:hiagwt}) in the late nineties.   

Yet another tessel in this mosaic is provided by moduli spaces of Admissible Covers, a smooth compactification of the classical Hurwitz schemes, parameterizing ramified covers of curves with prescribed numerical invariants and ramification data (degree of the cover, genus of the source and of the target curves, and ramification profile over all branch points). The natural forgetful map to the moduli space of curves (remembering the source) allows one to define Hodge-type integrals
on moduli spaces of admissible covers; such integrals were named Hurwitz-Hodge integrals by Jim Bryan, Tom Graber and Rahul Pandharipande. In their works (\cite{bgp:crc} and\cite{bg:prep}), inspired by Kevin Costello's \cite{c:thesis}, they pursue a systematic approach to the orbifold Gromov-Witten theory of Gorenstein stacks via Hurwitz-Hodge integrals. In \cite{bgp:crc}, through the explicit evaluation of $\lambda_{g-1}$ on moduli spaces of admissible $\mathbb{Z}_3$-covers, they describe the Gromov-Witten potential of the quotient stack $[\mathbb{C}^2/\mathbb{Z}_3]$ and show that the Crepant Resolution Conjecture is verified in that case.
 
We began studying Hurwitz-Hodge integrals in \cite{r:adm} and \cite{r:tqft}. In this paper we make some progress in the understanding of Hurwitz-Hodge integrals with descendants, exhibiting some surprisingly nice combinatorial structure. In \cite{bct:gg-1}, we use this structure to give a purely Hurwitz-theoretic proof of some classical computations of Faber-Pandharipande and Loojenga (see the ``Important Remark" in Section \ref{results} for a short discussion of this application).

\subsection{The Theorems} \label{results}
We fix once and for all the degree $d$ of the maps we consider.
For any positive integer $i$, define the generating function $\mathcal{D}_i(u)$ as follows:
$$
D_i^g := \int_{\mathcal{A}_g} \lambda_g\lambda_{g-i}\psi^{i-1},
$$
\begin{eqnarray}
\mathcal{D}_i(u) :=\sum_{g\geq i} D_i^g \frac{u^{2g}}{2g!}.	
\end{eqnarray}
\begin{theorem}
\label{d}
$$
\begin{imp}
\displaystyle{
\mathcal{D}_i(u)=\frac{d^{i-1}}{i!}\mathcal{D}_1^i(u).
}
\end{imp}
$$
\end{theorem}
Theorem \ref{d} suggests the following definition, that  fits nicely with the forthcoming computations.
 \begin{defi}
 \label{0}
 $$\mathcal{D}_0:= \frac{1}{d}.$$
\end{defi}
 
\textbf{Important Remark:} In this paper we consider the generating function $\mathcal{D}_1(u)$ as an initial condition. Such generating function was explicitly computed by Bryan-Pandaripande in \cite{bp:tlgwtoc}, who encoded in generating function language previous computations by Faber and Pandharipande  for $d=2$ (\cite{fp:lsahiittr}), and by Looijenga (\cite{l:ottr}) for all other degrees:
\begin{eqnarray}
\label{fapalo}
	\mathcal{D}_1(u)=\ln\left(\frac{d\sin\left(\frac{u}{2}\right)}{\sin\left(\frac{du}{2}\right)}\right).
\end{eqnarray}
 In \cite{bct:gg-1} with Aaron Bertram and Gueorgui Todorov, we provide a new proof of (\ref{fapalo}) relying only on Hurwitz theory and essentially exploiting the structure exhibited in Theorem \ref{d}.

Combining the result of  Theorem \ref{d}   with formula (\ref{fapalo}), it is straightforward to obtain the following:

\textbf{Corollary:}
\label{v}
$$
\begin{imp}
\displaystyle{
\sum_{g=0}^\infty \mathcal{T}_g(u) = \frac{1}{d}\ \mbox{e}^{d\mathcal{D}_1(u)} =\frac{d^{d-1}\sin^d\left(\frac{u}{2}\right)}{\sin^d\left(\frac{du}{2}\right)}.
}
\end{imp}
$$
 
\vspace{1cm}
For any non-negative integer $i$, define:
$$
V_i^g:= \int_{\mathcal{A}_g}\lambda_g\lambda_{i}\psi^{g-i-1}.
$$ 

Consistently with Definition \ref{0}, we define the unstable term:
$$
V_0^0:=\frac{1}{d}.
$$

Finally, the generating function:
\begin{eqnarray}
\mathcal{V}_i(u):=\sum_{g\geq i} V_i^g \frac{u^{2g}}{2g!}	
\end{eqnarray}

We  describe the generating function $\mathcal{V}_i(u)$ in terms of the initial conditions $V_i^{i+1}$ (that are the coefficients of the generating function $\mathcal{D}_1(u)$ defined above). Before we do so, we estabilish some notation.

\textbf{Notation:}
\begin{enumerate}
	\item To avoid carrying around double indices, we rename $V_i^{i+1}= \Upsilon_i$.
	\item \textbf{Partition notation:} For $\eta=n_1^{m_1}\ldots n_r^{m_r}$ a partition of the integer $i$ ($\eta \vdash i$), we denote 
\begin{enumerate}
	\item $|\eta|= i = \sum m_kn_k$.
	\item $\ell(\eta) = \sum m_k$.
	\item $Aut(\eta) =\prod m_k!$.
	\item $\Upsilon^\eta = \Upsilon_{n_1}^{m_1}\cdot \ldots \cdot \Upsilon_{n_r}^{m_r}$.
\end{enumerate}
\end{enumerate}

Finally we are ready to state:
\begin{theorem}
\label{v}
$$
\begin{imp}
\displaystyle{
\mathcal{V}_i (u) = u^{2i}\mbox{e}^{d \Upsilon_0 u^2}\sum_{\eta \vdash i}u^{2\ell(\eta)} d^{\ell(\eta)-1}\frac{\Upsilon^{\eta}}{Aut(\eta)}.
}
\end{imp}
$$
\end{theorem}

\subsection{Plan of the Paper}

The paper is organized as follows:
\begin{itemize}
	\item Sections \ref{ac} and \ref{loc} are meant to provide background and establish notation. In Section \ref{ac} we quickly review the basic facts about moduli spaces of admissible covers that we subsequently use. Section \ref{loc} is a minimalistic presentation of localization, mainly focused on showing how it is used in this particular application.
	\item Section \ref{thm1} presents the proof of Theorem \ref{d}.
	\item Section \ref{thm2} gives the proof of Theorem \ref{v}. Since in this case the combinatorial aspect or the proof is pretty interesting and somewhat sophisticated, we have chosen to split the proof into two parts: ``Geometry" and ``Combinatorics".
\end{itemize}
\section{Admissible Covers}\label{ac}

\subsection{Basic Definitions and Notation} \label{adm}
Moduli spaces of admissible covers are a ``natural'' compactification of
the Hurwitz schemes, parameterizing ramified covers of smooth Riemann
Surfaces.  Let $ (X,p _1,\cdots,p_r) $ be an $ r $-pointed
nodal curve of genus $ g $.
\begin{defi}
    An \textbf{admissible cover} $ \pi:E\longrightarrow X $ of degree d
    is a finite morphism satisfying the following:
    \begin{enumerate}
        \item
            E is a nodal curve.
        \item
            Every node of E maps to a node of X.
        \item
            The restriction of $ \pi:E\longrightarrow X $ to $ X\smallsetminus
            (p_1,\cdots,p_r) $ is \'{e}tale of constant degree d.
        \item
            Nodes can be smoothed. This means: given an admissible cover
$\pi:E\rightarrow X$, and a node of $E$, we can find a family of admissible 
covers $\pi':E'\rightarrow X'$ such that:
\begin{itemize}
\item $\pi:E\rightarrow X$ is the central fiber of the family;
\item locally in analytic coordinates, X', E' and $ \pi'
            $ are described as follows, for     some positive integer $n$ not larger than $d$:
            $$
                \begin{array}{cl}
                    E:  & e_1e_2=a, \\
                    X:  & x_1x_2=a^n ,\\
                    \pi:  & x_1= e_1^n,\ x_2= e_2^n.
                \end{array}
            $$
\end{itemize}

    \end{enumerate}
\end{defi}

\textbf{Notation:} In this paper we are concerned with moduli spaces of admissible covers of degree $d$ of an unparameterized $\proj$, satisfying the following ramification data:
\begin{itemize}
	\item over the first $2g$ marked points we have simple ramification ($t$ for transposition - the monodromy type at each of these marked points).
		\item we require the profile of the cover over the last two marked points to consist of only one point. This corresponds to full ramification, or monodromy type given by a $d$-cycle (hence $(d)$).
\end{itemize}
  The notation we adopt is maybe a little cumbersome, but it has the advantage of containing in an unequivocal fashion all the combinatorial information we are working with:
$$
\A{g}{0}.
$$
This moduli space is a smooth DM stack of dimension $2g-1$.

\textbf{Remark:} In the language of \cite{acv:ac}, we are selecting the connected components of the stack of twisted stable maps $\mathcal{K}_{0,2g+2}(BS_d,0)$ satisfying the above ramification conditions.

\subsection{The boundary}\label{nodal}

Admissible covers of a nodal curve can be combinatorially described in terms of 
admissible covers of the irreducible components of the curve. This is extremely 
useful because it opens the way to the use of degeneration techniques and induction.
Crucial are the following identities (\cite{l:adffgwi}), that take place in the Chow ring
with rational coefficients.

Let \{A,B\} be a two set partition of the set of $2g+2$ marks on the base curve, and denote by 
$$
D(A\mid B) \subset \A{g}{0}
$$
the divisor corresponding to the base curve splitting into a nodal rational curve with (at least) one node, the marks in set $A$ arranging themselves on one component, those in $B$ on the other.

Then: 
\begin{eqnarray}
        \displaystyle{[D(A\mid B)]=
        \sum_{\eta \vdash d}\mathfrak{z}(\eta)[\overline{Adm}_{g_1\stackrel{d}{\rightarrow}0,
        (A, {\eta})}] \times [ \overline{Adm}_{g_2\stackrel{d}{\rightarrow}0,
        (B,{\eta})}],} 
\label{rnodal}
 \end{eqnarray}
where:
\begin{itemize}
        \item $\eta = ((\eta^1)^{m_1},\ldots,(\eta^k)^{m_k})$ runs over all partitions
 of $d$;
        \item  the combinatorial factor
\begin{eqnarray}
    \mathfrak{z}(\eta):=\prod m_i!(\eta^i)^{m_i}
    \label{zeta}.
\end{eqnarray}
 is the order of the centralizer in $S_d$ of any group element in the conjugacy class of $\eta$;
        \item $g_1$ and $g_2$ are determined by the Riemann-Hurwitz formula and they satisfy:
        $$g_1+g_2+\ell(\eta)-1=g.$$
\end{itemize}

\subsection{Tautological Classes} \label{taut}
Moduli spaces of admissible covers admit two natural forgetful maps as in the following diagram:
$$
\begin{array}{ccc}
\A{g}{0} & \stackrel{s}{\rightarrow} & \overline{M}_g \\ & & \\
\mbox{\footnotesize $t$}\downarrow & & \\ & & \\
\overline{M}_{0,2g+2}.
\end{array}
$$

The vertical map remembers the base of the cover as a genus $0$ curve marked by the branch points of the cover.
On $ \overline{M}_{0,2g+2}$ the tautological class $\psi_i\in A^1(\overline{M}_{0,2g+2}) $ is the first Chern class of the cotangent line bundle $\mathbb{L}_i$ (\cite{k:pc}).
\begin{defi}
    The tautological class $ \psi_i^{Adm} \in A^1(\A{g}{0}) $ is defined to be the
    pull-back of the analogous class via the map $ t $:
    $$
        \psi_i^{Adm}:=t^\ast(\psi_i).
    $$
\end{defi}

The map $s$ forgets the cover map and only remembers the source curve.
The tautological class $ \lambda_i \in A^i(\overline{M}_{g}) $ is
defined to be the $ i $-th Chern class of the Hodge bundle $ \mathbb{E}$. 
\begin{defi}
    The tautological class $ \lambda_i^{Adm} \in A^i(\overline{Adm}_{h\stackrel
    {d}{\rightarrow}0, (\eta_1,\cdots,\eta_{n})}) $ is defined to be the pull-back of the corresponding 
    class on the moduli space of curves:
    $$
        \lambda_i^{Adm}:=s^\ast (\lambda_i).
    $$
\end{defi}

We immediately drop the superscripts ``Adm'' since we are only dealing with admissible cover classes.

\subsubsection*{Aside: on the Hodge bundle and its useful properties.}

The Hodge bundle $\mathbb{E}$ \footnote{we will add a subscript and denote it $\mathbb{E}_g$ when we want to emphasize the genus.} is a natural rank $g$ bundle on the moduli space of curves $\overline{M}_{g}$. It is defined to be the pushforward via the universal family map of the relative dualizng sheaf:
$$
\mathbb{E} =\pi_\ast {\omega}_\pi.
$$
Over a smooth curve, we can interpret the fiber of $\mathbb{E}$ as $h^0(C,K_C)$, the space of global holomorphc differentials on $C$.

The following properties of $\mathbb{E}$ are important and useful(\cite{m:taegotmsoc}):
\begin{description}
	\item[Mumford Relation:] the total Chern class of the sum of the Hodge bundle with its dual is trivial:
	\begin{eqnarray}\label{mumf}c(\mathbb{E}\oplus \mathbb{E}^\ast)=1.\end{eqnarray}
	\item[Top Chern Class:] an immediate consequence of the previous point is that for all $g>0$
	\begin{eqnarray}\label{gsquared}\lambda_g^2=0.\end{eqnarray}
	\item[Separating nodes:] denote by $\Delta_{i,g-i}$ the divisor in $\overline{M}_{g}$ parameterizing nodal curves $C=C_1\bigcup_p C_2$, with $C_1$ a curve of genus $i$, $C_2$ of genus $g-i$. Then the restriction of the Hodge bundle splits as the sum of the bundles corresponding to the two components:
	\begin{eqnarray}\label{disconnect}\mathbb{E}\mid_{\Delta_{i,g-i}}=\mathbb{E}_i \oplus \mathbb{E}_{g-i}.\end{eqnarray}
		\item[Non-separating nodes:] denote by $\Delta_0$ the divisor in $\overline{M}_{g}$ parameterizing nodal curves obtained by attaching two points of a genus $g-1$ curve. Then the restriction of the Hodge bundle splits a trivial factor:
	\begin{eqnarray}\label{selfnode}\mathbb{E}\mid_{\Delta_0}=\mathbb{E}_{g-1} \oplus \mathcal{O}.\end{eqnarray}
	\end{description}

\subsection{Admissible Covers of a Parametrized $\proj$}

These are a  variation of the previous moduli spaces: the objects are parametrized are the same, but the equivalence
relation is stricter:  we consider two covers $ E_1\rightarrow\proj $, $
E_2\rightarrow\proj $ equivalent if there is an isomorphism $ \varphi:E_1\rightarrow
E_2 $ that makes the natural triangle commute.  In other words, we are
not allowed to act on the base with an automorphism of $ \proj $.

We denote by
$$
\A{g}{\proj}
$$ 
the stack of admissible covers of genus $g$ and degree $d$ of a parametrized projective line, with specified ramification data. This is a smooth DM stack of dimension $2g+2$.

\textbf{Remark:} in the language of \cite{acv:ac}, we are selecting the connected components satisfying the appropriate ramification conditions in the stack of twisted stable maps $\mathcal{K}_{0,2g+2}([\proj/S_d],d!)$; the symmetric group acts trivially on $\proj$.

When the base curve degenerates to a nodal rational curve, one component will remain parametrized, whereas the sprouted twigs will be  unparameterized genus $0$ curves. Therefore ordinary genus $0$ admissible cover spaces naturally appear in the combinatorial description of the boundary strata of parametrized admissible cover spaces.

\section{Localization}\label{loc}

Consider the  one-dimensional algebraic torus $\mathbb{C}^\ast$, and recall that the $\mathbb{C}^\ast$-equivariant Chow ring of a point is a polynomial ring in one variable:
$$A^\ast_{\Cstar}(\{pt\},\mathbb{C})= \mathbb{C}[\hbar]. $$
Let $\Cstar$ act on a smooth, proper stack $X$, denote by $i_k:F_k\hookrightarrow X$ the irreducible components of the fixed locus for this action and by $N_{F_k}$ their normal bundles. The natural map:
$$
\begin{array}{ccc}
A^\ast_{\Cstar}(X) \otimes \mathbb{C}(\hbar) & \rightarrow & \sum_{k}{A^\ast_{\Cstar}}(F_k) \otimes \mathbb{C}(\hbar)\\
                                             &             &                                        \\
\alpha                                       & \mapsto     &\displaystyle{\frac{i_k^\ast\alpha}{c_{top}(N_{F_k})}}.
\end{array}
$$
is an isomorphism. Pushing forward equivariantly to the class of a point, we obtain the Atiyah-Bott integration formula:
$$\int_{[X]}\alpha = \sum_k \int_{[F_k]} \frac{i_k^\ast\alpha}{c_{top}(N_{F_k})}.$$

\subsection{Our Set-up}
Let $\mathbb{C}^\ast$ act on a two-dimensional vector space $V$ via:
$$t\cdot(z_0,z_1)=(tz_0,z_1).$$
This action descends on $\proj$, with fixed points $0=(1:0)$ and $\infty=(0:1)$. An equivariant lifting of $\Cstar$ to a line bundle $L$ over $\proj$ is uniquely determined by its  weights $\{L_0,L_\infty\}$ over the fixed points.

The canonical lifting of $\Cstar$ to the tangent bundle of $\proj$ has weights $\{1,-1\}$.

The action on $\proj$ induces an action on the moduli spaces of admissible covers to a parametrized $\proj$ simply by post composing the cover map with the automorphism of $\proj$ defined by $t$. 

The fixed loci for the induced action on the moduli space consist of admissible covers such that anything ``interesting'' (ramification, nodes, marked points) happens over $0$ or $\infty$, or on ``non special'' twigs that attach to the main $\proj$ at  $0$ or $\infty$.

\section{Theorem 1}
\label{thm1}

We compute the generating functions $\mathcal{D}_i(u)$  by evaluating via localization the following auxiliary integrals:
\begin{eqnarray}
	I^g_i:= \int_{\mathcal{A}_g} \lambda_g\lambda_{g-i}\ ev_1^{\ast}(0)\  ev_{2g+1}^{\ast} (0)\ ev_{2g+2}^{\ast}(\infty),
	\label{aux}
\end{eqnarray}

where:
\begin{itemize}
	\item $\mathcal{A}_g =\A{g}{\proj}$;
	\item  $0\leq i \leq g$.
\end{itemize}
 \textbf{Remark:} The integrals $I^g_0$ vanish for $g> 0$ because of Mumford's relation (\ref{gsquared}). It is straightforward to compute 
 $$I^0_0= \frac{1}{d}.$$
 This computation is consistent with our choice of defining  $\mathcal{D}_0=1/d$ (Definition \ref{0}).
 
\begin{figure}[htbp]
	\begin{center}
$$	
\begin{array}{ccc}
	F_g = \ \ 	\includegraphics[width=0.32\textwidth]{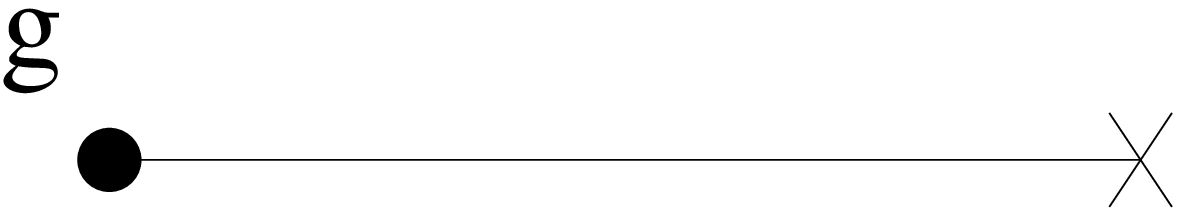} &
				\hspace{.7cm}																							&		
	F_{g_1g_2}=\  \	\includegraphics[width=0.35\textwidth]{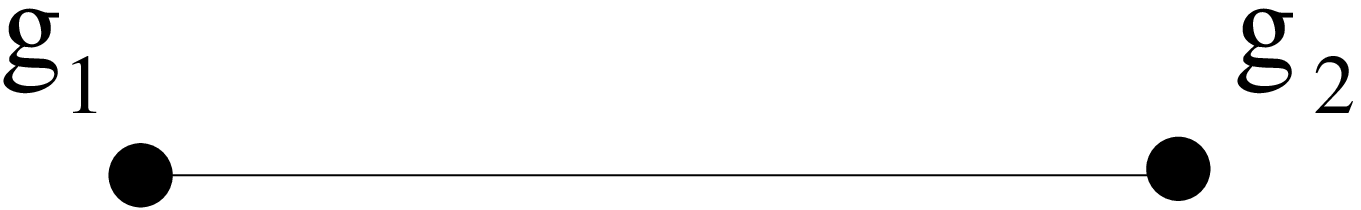}
\end{array}
$$
	\end{center}
	\caption{The possibly contributing fixed loci.}
	\label{fix}
\end{figure}

When we evaluate $I^g_i$ via localization, the number of possibly contributing fixed loci (identified here with their localization graphs) is reduced because of the following considerations:
\begin{enumerate} 
	\item The full ramification condition imposed both at $0$ and $\infty$ forces the preimages of $0$ and $\infty$ to be connected. 
	\item The presence of the class $\lambda_g$ in the integrand forces no loops in the localization graph (from relation (\ref{selfnode})).
	\item The additional simple transposition  ``sent" to $0$ rules out the fixed locus where the preimage of $0$ is just a single point.
\end{enumerate}
The possibly contributing fixed loci consist of boundary strata parameterizing a single sphere, fully ramified over $0$ and $\infty$, mapping with degree $d$ to the main $\proj$. A rational tail must sprout from $0$, and be covered by a curve of genus $g_1$. If $g_1<g$, then a rational tail sprouts from $\infty$ as well, covered by a curve of genus $g_2=g-g_1$.
The situation is illustrated in Figure \ref{fix}.
  
\subsubsection*{The case i=1}
\label{case1}
In this case it is easy to see, by pure dimension reasons, that the only contributing fixed locus is the only codimension three locus, that we have denoted in Figure \ref{fix} by $F_g$.
Noting that $$F_g \cong \A{g}{0}\ \ ,$$ the explicit evaluation of the integral yields:
\begin{eqnarray}
I^g_1 =\int_{F_g}\frac{-\hbar^3}{\hbar (\hbar-\psi) (-\hbar)}\lambda_g\lambda_{g-1}= D_1^g.
\label{gg-1}
\end{eqnarray}

\subsubsection*{The case $i\geq 2$}
\label{case2}
   For $i \geq 2$ the integral $I(g,i)$ vanishes for dimension reasons: we are integrating a $(2g-i+3)$-dimensional class on a $(2g+2)$-dimensional space. Localization produces inductive relations between our generating functions.
The fixed locus $F_g$ behaves differently from all the $F_{g_1g_2}$. For this reason we analyze their contributions separately.
\begin{description}
	\item[$F_g$: ] 
	$$
	\int_{F_g}\frac{-\hbar^3}{\hbar (\hbar-\psi) (-\hbar)}\lambda_g\lambda_{g-i}=\frac{1}{\hbar^{i-1}} \int_{F_g} \lambda_g\lambda_{g-i}\psi^{i-1} =D_i^g.
	$$
	\item[$F_{g_1g_2}$: ] this fixed locus is isomorphic to a product of spaces, with multiplicities:
	$$
	\hspace{-1.3cm}F_{g_1g_2}\cong d {{2g-1}\choose{2g_1-1}} \A{g_1}{0} \times \A{g_2}{0}.
	$$
	
\begin{itemize}
	\item The factor of $d$ is a combination of the $d^2$ coming from gluing two nodes ((\ref{rnodal}),(\ref{zeta})) and the $1/d$ automorphism contribution of the fully ramified $\proj$ lying over the main $\proj$;
	\item The combinatorial factor keeps track of all possible combinations of $2g_1-1$ marks ``going to $0$" and $2g_2$ ``going to $\infty$".
\end{itemize}
With this in mind:
\begin{eqnarray}
\hspace{-1.3cm}	\int_{F_{g_1,g_2}}\frac{-\hbar^3}{\hbar (\hbar-\psi_0) \hbar(\hbar+\psi_\infty)}(\lambda_{g_1}\lambda_{g_2})(\lambda_{g_1-i}\lambda_{g_2} +\lambda_{g_1-i+1}\lambda_{g_2-1} +\ldots +\lambda_{g_1}\lambda_{g_2-i}) \nonumber\\
=\frac{1}{\hbar^{i-1}}d{{2g-1}\choose{2g_1-1}} \sum_{k=1}^i (-1)^kD^{g_1}_{i-k}D^{g_2}_{k}. 
\label{g1g2}
\end{eqnarray}
\end{description}
Adding the contributions of all fixed loci, and remembering that the integral $I_i^g=0$ we obtain the relation\footnote{Note that a priori $g_1,g_2 \not =0$, but we can omit this condition because $D_k^0$ is $0$.}:
\begin{eqnarray}
D_i^g = -d\sum_{g_1+g_2=g}{{2g-1}\choose{2g_1-1}}\sum_{k=1}^i (-1)^kD^{g_1}_{i-k}D^{g_2}_{k}.
\label{relation}
\end{eqnarray}
\textbf{Remark:} Relation \ref{relation} determines $\mathcal{D}_i (u)$, in terms of all $\mathcal{D}_j (u)$ with $j<i$. 
After interchanging the summation order, and dividing out by $(2g-1)!$, (\ref{relation}) expresses the following identity of generating functions:
\begin{eqnarray}
\label{gfrelation}
\frac{d}{du}\mathcal{D}_i = -d \sum_{k=1}^i (-1)^k\left(\frac{d}{du}\mathcal{D}_{i-k}\right)\mathcal{D}_{k}.
\end{eqnarray}
Recalling that we defined $\mathcal{D}_0=1/d$, we obtain an even simpler expression for (\ref{gfrelation}):
\begin{eqnarray}
\label{gfrelationnew}
	\begin{imp}
	\displaystyle{
\sum_{k=0}^i (-1)^k\left(\frac{d}{du}\mathcal{D}_{i-k}\right)\mathcal{D}_{k}=0
	}
	\end{imp}
\end{eqnarray}
Theorem \ref{d} is  proved by plugging the predicted formula  in  (\ref{gfrelationnew}). The relation is immediately verified upon recognizing the elementary combinatorial identity:
$$
\sum_{k=0}^{i-1} (-1)^{k-1} {{i-1}\choose{k}}=0
$$

\section{Theorem 2}\label{thm2}
\subsection{Geometry}
\label{geom}
The strategy for the proof of Theorem \ref{v} is very similar. After all, we are studying the same information, only ``packaged" in a different way. With all notation as in (\ref{aux}), we define and localize the auxiliary integral:
\begin{eqnarray}
	J^g_i:= \int_{\mathcal{A}_g} \lambda_g\lambda_{i}\ ev_1^{\ast}(0)\  ev_{2g+1}^{\ast} (0)\ ev_{2g+2}^{\ast}(\infty).
\end{eqnarray} 
For $g>i+1$, this integral vanishes for dimension reasons. The discussion of the contributing fixed loci is completely analogous to Section \ref{thm1}, and the contributions are:
\begin{description}
	\item[$F_g$: ] 
	$$
	\int_{F_g}\frac{-\hbar^3}{\hbar (\hbar-\psi) (-\hbar)}\lambda_g\lambda_{i}=\frac{1}{\hbar^{i-1}} \int_{F_g} \lambda_g\lambda_{i}\psi^{i-1} =V_i^g.
	$$
	\item[$F_{g_1g_2}$: ]
\begin{eqnarray}
	\int_{F_{g_1,g_2}}\frac{-\hbar^3}{\hbar (\hbar-\psi_0) \hbar(\hbar+\psi_\infty)}(\lambda_{g_1}\lambda_{g_2})(\lambda_{i}1 +\lambda_{i-1}\lambda_{1} +\ldots +1\lambda_{i}) \nonumber\\
=\frac{1}{\hbar^{i-1}}d{{2g-1}\choose{2g_1-1}} \sum_{k=0}^i (-1)^{g_2-k}V^{g_1}_{k}V^{g_2}_{i-k}. 
\label{g1g2}
\end{eqnarray}
\end{description}
Adding the contributions over all fixed loci and recalling $J_i^g = 0$ (for $g>i+1$), we obtain:
\begin{eqnarray}
	\label{rel2}
	V_i^g+d\sum_{g_1+g_2=g}{{2g-1}\choose{2g_1-1}} \sum_{k=0}^i (-1)^{g_2-k}V^{g_1}_{k}V^{g_2}_{i-k}=0.
\end{eqnarray}

\textbf{Remark:} Note that the summation in relation (\ref{rel2}) occurs for $g_1$ and $g_2$ strictly positive integers. Hence (\ref{rel2}) determines the value of $V_i^g$ inductively in terms of:
\begin{itemize}
	\item the generating functions $\mathcal{V}_j$, with $j<i$;
	\item the initial condition $V_i^{i+1}$.
\end{itemize}

Relation (\ref{rel2}) assumes a particularly compact shape in generating function form thanks to the introduction of the unstable term $V_0^0=1/d$ in $\mathcal{V}_0$. After only a little bit of careful bookkeeping, we recognize (\ref{rel2}) to be the coefficient of $u^{2g-1}$ of the identity:
\begin{eqnarray}
\begin{imp}
\displaystyle{
	\sum_{k=0}^i (-1)^{i-k}\left(\frac{d}{du}\mathcal{V}_{k}(u)\right) \mathcal{V}_{i-k} (\im u) = \frac{(2i+2)}{d}V_i^{i+1}u^{2i+1}
}
\label{rel2good}
\end{imp}
\end{eqnarray}
\subsection{Combinatorics} 
\label{comb}

In this section we show that the generating functions
$$
\mathcal{V}_i (u) = u^{2i}\mbox{e}^{d \Upsilon_0 u^2}\sum_{\eta \vdash i}u^{2\ell(\eta)} d^{\ell(\eta)-1}\frac{\Upsilon^{\eta}}{Aut(\eta)}
$$
defined in page \pageref{v}, verify (\ref{rel2good}), thus concluding the proof of Theorem \ref{v}. 

Plugging these expressions in the left hand side of (\ref{rel2good}), we obtain:
\begin{eqnarray}
	u^{2i}\sum_{\eta \vdash i} \Upsilon^{\eta} \sum_{\eta_1\coprod \eta_2=\eta}\left[ 
	2\Upsilon_0 u^{2\ell(\eta)+1} d^{\ell(\eta)-1} (-1)^{\ell(\eta_2)}\frac{1}{Aut(\eta_1)}\frac{1}{Aut(\eta_2)}+ \right. \nonumber\\ 
\left.	+u^{2\ell(\eta)-1} d^{\ell(\eta)-2} (-1)^{\ell(\eta_2)}(2|\eta_1|+2\ell(\eta_1))\frac{1}{Aut(\eta_1)}\frac{1}{Aut(\eta_2)}
	\right]
\end{eqnarray}
The theorem now follows immediately from these two purely combinatorial statements.

\begin{lemma}
\label{comb1}
\begin{eqnarray}
	\sum_{\eta_1\coprod \eta_2=\eta}(-1)^{\ell(\eta_2)}\frac{1}{Aut(\eta_1)}\frac{1}{Aut(\eta_2)} =0.
\label{comb1f}
\end{eqnarray}
\end{lemma}
\begin{lemma}
\label{comb2}
\begin{eqnarray}\label{comb2f}
	\begin{array}{ccc}
		\displaystyle{\sum_{\eta_1\coprod \eta_2=\eta}(-1)^{\ell(\eta_2)}(|\eta_1|+\ell(\eta_1))\frac{1}{Aut(\eta_1)}\frac{1}{Aut(\eta_2)}} & = &
	\left\{
		\begin{array}{lc}
		(i+1) & \mathrm{if}\ \eta= (i), \\
		          &                     \\
		0         & \mathrm{else.}
		\end{array}
	\right.
	\end{array}
\end{eqnarray}
\end{lemma}
\textsc{Proof of Lemma \ref{comb1}:} Proving this statement amounts to recognizing expression (\ref{comb1f}) as a product of binomial coefficients. If $\eta= n_1^{m_1}\ldots n_r^{m_r}$, then
$$
\prod_{j=1}^r (X_j-1)^{m_j}=\prod_{j=1}^r\left(\sum_{k_j=0}^{m_j} {{m_j}\choose{k_j}} (-1)^{m_j-k_j} X_j^{k_j}\right)=
$$
$$
=Aut(\eta)\sum_{\eta_1\coprod \eta_2=\eta}(-1)^{\ell(\eta_2)}\frac{\prod X_j^{k_j}}{Aut(\eta_1)}\frac{1}{Aut(\eta_2)},
$$
where $\eta_1= n_1^{k_1}\ldots n_r^{k_r}$. Now Lemma (\ref{comb1}) follows by setting all of the $X_j=1$.

\vspace{1cm}
\textsc{Proof of Lemma \ref{comb2}:} We observe first of all that for $\eta= (i)$, the statement is easily checked. 
Now let $\eta =n_1^{m_1}\ldots n_r^{m_r}$. Denote 
$$
\tilde{\eta} =n_1^{m_1}\ldots n_{r-1}^{m_{r-1}}.
$$
Any subdivision $\eta=\eta_1 \coprod \eta_2$ induces a subdivision $\tilde{\eta}=\tilde{\eta}_1 \coprod \tilde{\eta}_2$ simply by forgetting the $n_r$ parts of the partition $\eta$. We now group the the terms of (\ref{comb2f}) according to this induced subdivisions:
$$
\mbox{LHS of}\ (\ref{comb2f})=
$$
$$
=\sum_{\tilde{\eta}_1\coprod \tilde{\eta}_2=\tilde{\eta}}\sum_{k=0}^{m_r}(-1)^{\ell(\tilde{\eta}_2)+m_r -k}(|\tilde{\eta}_1|+ kn_r +\ell(\tilde{\eta}_1)+k)\frac{1}{Aut(\tilde{\eta}_1)}\frac{1}{k!}\frac{1}{Aut(\tilde{\eta}_2)}\frac{1}{(m_r-k)!}=
$$
$$
=
\left
(\sum_{\tilde{\eta}_1\coprod \tilde{\eta}_2=\tilde{\eta}}(-1)^{\ell(\tilde{\eta}_2)}	(|\tilde{\eta}_1|+\ell(\tilde{\eta}_1))\frac{1}{Aut(\tilde{\eta}_1)}\frac{1}{Aut(\tilde{\eta}_2)} 
\right)
\left(
\sum_{k=0}^{m_r}(-1)^{m_r-k}  \frac{1}{k!}\frac{1}{(m_r-k)!}
\right)+
$$
$$
+(n_r+1)\left
(\sum_{\tilde{\eta}_1\coprod \tilde{\eta}_2=\tilde{\eta}}(-1)^{\ell(\tilde{\eta}_2)}	\frac{1}{Aut(\tilde{\eta}_1)}\frac{1}{Aut(\tilde{\eta}_2)} 
\right)
\left(
\sum_{k=0}^{m_r}(-1)^{m_r-k}k  \frac{1}{k!}\frac{1}{(m_r-k)!}
\right)
$$

Let us  observe the summations in this expression:
\begin{description}
	\item[$\displaystyle{\sum_{k=0}^{m_r}(-1)^{m_r-k}  \frac{1}{k!}\frac{1}{(m_r-k)!}}$:] this term is clearly always $0$;  we recognize (up to sign and a global factor of  $m_r!$) the binomial expansion of $(1-1)^{m_r}$.
\item[$\displaystyle{
\sum_{\tilde{\eta}_1\coprod \tilde{\eta}_2=\tilde{\eta}}(-1)^{\ell(\tilde{\eta}_2)}	\frac{1}{Aut(\tilde{\eta}_1)}\frac{1}{Aut(\tilde{\eta}_2)} 
}$:] this vanishing is precisely the statement of Lemma \ref{comb1}.
\end{description} 

This concludes the proof of the Lemma and of Theorem \ref{v}.

\textbf{Remark:} It is interesting to observe that also the other two factors in the above expression vanish ``almost" always:
\begin{description}
\item[$\displaystyle{\sum_{\tilde{\eta}_1\coprod \tilde{\eta}_2=\tilde{\eta}}(-1)^{\ell(\tilde{\eta}_2)}	(|\tilde{\eta}_1|+\ell(\tilde{\eta}_1))\frac{1}{Aut(\tilde{\eta}_1)}\frac{1}{Aut(\tilde{\eta}_2)} 
}$:] this term vanishes unless $\ell(\tilde{\eta})=1$. This is precisely the statement of Lemma \ref{comb2} for the partition $\tilde{\eta}$, and it can be shown by induction on the size of the partition.
\item[$\displaystyle{\sum_{k=0}^{m_r} (-1)^{m_r-k}k \frac{1}{k!}\frac{1}{(m_r-k)!}}$:]if we multiply this summation by $m_r!$, we recognize the binomial expansion of 
	$$
	\frac{d}{dx}(x-1)^{m_r},
	$$
	evaluated at $x=1$. The vanishing follows for $m_r>1$.
	
\end{description}


\addcontentsline{toc}{section}{Bibliography}
\bibliographystyle{alpha}
\bibliography{biblio}

\begin{thebibliography}{ELSV01}

\bibitem[ACV01]{acv:ac}
Dan Abramovich, Alessio Corti, and Angelo Vistoli.
\newblock Twisted bundles and admissible covers.
\newblock {\em Comm in Algebra}, 31(8):3547--3618, 2001.

\bibitem[BCT06]{bct:gg-1}
Aaron Bertram, Renzo Cavalieri, and Gueorgui Todorov.
\newblock Evaluating tautological classes using only {H}urwitz numbers.
\newblock Preprint:on the ArXiv soon., 2006.

\bibitem[BG05]{bg:prep}
Jim Bryan and Tom Graber.
\newblock The crepant resolution conjecture.
\newblock In preparation, 2005.

\bibitem[BGP05]{bgp:crc}
Jim Bryan, Tom Graber, and Rahul Pandharipande.
\newblock The orbifold quantum cohomology of $\mathbb{C}^2/\mathbb{Z}_3$ and
  {Hurwitz-Hodge} integrals.
\newblock Preprint:math.AG/0510335, 2005.

\bibitem[BP04]{bp:tlgwtoc}
Jim Bryan and Rahul Pandharipande.
\newblock The local {Gromov-Witten} theory of curves.
\newblock Preprint: math.AG/0411037, 2004.

\bibitem[Cav04]{r:adm}
Renzo Cavalieri.
\newblock Hodge-type integrals on moduli spaces of admissible covers.
\newblock Preprint: math.AG/0411500, 2004.

\bibitem[Cav05]{r:tqft}
Renzo Cavalieri.
\newblock A {TQFT} for intersection numbers on moduli spaces of admissible
  covers.
\newblock Preprint: mathAG/0512225, 2005.

\bibitem[Cos03]{c:thesis}
Kevin Costello.
\newblock Higher-genus {G}romov-{W}itten invariants as genus 0 invariants of
  symmetric products.
\newblock Preprint:mathAG/0303387, 2003.

\bibitem[ELSV99]{elsv:ohnahi}
Torsten Ekedahl, Sergei Lando, Michael Shapiro, and Alek Vainshtein.
\newblock On {Hurwitz} numbers and {Hodge} integrals.
\newblock {\em C.R. Acad.Sci.Paris Ser.I Math.}, 328:1175--1180, 1999.

\bibitem[ELSV01]{elsv:hnaiomsoc}
Torsten Ekedahl, Sergei Lando, Michael Shapiro, and Alek Vainshtein.
\newblock {Hurwitz} numbers and intersections on muduli spaces of curves.
\newblock {\em Invent. Math.}, 146:297--327, 2001.

\bibitem[Fab99]{f:algo}
Carel Faber.
\newblock Algorithms for computing intersection numbers on moduli spaces of
  curves, with an application to the class of the locus of jacobians.
\newblock {\em New trends in algebraic geometry (Warwick 1996),London Math.
  Soc. Lecture Note Ser.}, 264:93--109, 1999.

\bibitem[FP00a]{fp:lsahiittr}
C.~Faber and R.~Pandharipande.
\newblock Logarithmic series and {H}odge integrals in the tautological ring.
\newblock {\em Michigan Math. J.}, 48:215--252, 2000.
\newblock With an appendix by Don Zagier, Dedicated to William Fulton on the
  occasion of his 60th birthday.

\bibitem[FP00b]{fp:hiagwt}
Carel Faber and Rahul Pandharipande.
\newblock {Hodge} integrals and {Gromov-Witten} theory.
\newblock {\em Invent. Math.}, 139(1):173--199, 2000.

\bibitem[GJV03]{gjv:ttgodhn}
I.~Goulden, D.~M. Jackson, and Ravi Vakil.
\newblock Towards the geometry of double {Hurwitz} numbers.
\newblock Preprint: math.AG/0309440v1, 2003.

\bibitem[GJV06]{gjv:spolgc}
Ian Goulden, David Jackson, and Ravi Vakil.
\newblock A short proof of the $\lambda_g$-conjecture without {G}romov-{W}itten
  theory: Hurwitz theory and the moduli of curves.
\newblock Preprint:mathAG/0604297, 2006.

\bibitem[GV01]{gv:taut}
Tom Graber and Ravi Vakil.
\newblock On the tautological ring of {$\overline{M}_{g,n}$}.
\newblock {\em Turkish J. Math.}, 25(1):237--243, 2001.

\bibitem[GV03a]{gv:hnavl}
Tom Graber and Ravi Vakil.
\newblock {Hodge} integrals, {Hurwitz} numbers, and virtual localization.
\newblock {\em Compositio Math.}, 135:25--36, 2003.

\bibitem[GV03b]{gv:rvl}
Tom Graber and Ravi Vakil.
\newblock Relative virtual localization and vanishing of tautological classes
  on moduli spaces of curves.
\newblock Preprint: math.AG/0309227, 2003.

\bibitem[Koc01]{k:pc}
Joachim Kock.
\newblock Notes on psi classes.
\newblock Notes, 2001.

\bibitem[Li02]{l:adffgwi}
Jun Li.
\newblock A degeneration formula of {GW}-invariants.
\newblock {\em J. Differential Geom.}, 60(2):199--293, 2002.

\bibitem[Loo95]{l:ottr}
Eduard Looijenga.
\newblock On the tautological ring of {${ M}_g$}.
\newblock {\em Invent. Math.}, 121(2):411--419, 1995.

\bibitem[Mum83]{m:taegotmsoc}
David Mumford.
\newblock Toward an enumerative geometry of the moduli space of curves.
\newblock {\em Arithmetic and Geometry}, II(36):271--326, 1983.

\end{thebibliography}

\end{document}